\documentclass[11pt,reqno]{amsart}
\usepackage{amsfonts}

\usepackage{amsaddr}
\usepackage{graphicx}
\usepackage{amscd}
\usepackage{amsthm}
\usepackage{amstext}
\usepackage{amsmath}
\usepackage{float}
\usepackage{amssymb,latexsym}
\usepackage{epsfig}
\usepackage{mathrsfs}
\usepackage{subfigure}
\usepackage{subfig}
\usepackage{cite}
\usepackage{multicol}
\usepackage[margin=2.5cm]{geometry}
\DeclareGraphicsExtensions{.png,.pdf,.jpg,.gif}

\newtheorem{teo}{Theorem}[section]

\theoremstyle{definition}

\newtheorem{remark}{Remark}[section]

\theoremstyle{plain}

\numberwithin{equation}{section}

\begin{document}

\title[Periodic solutions and a small parameter] {Bifurcation and periodic solutions to neuroscience models with a small parameter}

\author[Jos\'e Oyarce]{Jos\'e Oyarce \textsuperscript{$\MakeLowercase{a}$}}
\email{jooyarce@egresados.ubiobio.cl}

\address{$^a$Departamento de Matem\'atica, Facultad de Ciencias, Universidad del B\'io-B\'io, Casilla 5-C, Concepci\'on, VIII-Regi\'on, Chile}


\begin{abstract}
The existence of periodic solutions is proven for some neuroscience models with a small parameter. Moreover, the stability of such solutions is investigated, as well. The results are based on a theoretical research dealing with the functional differential equation with parameters
$$
\dot{x}(t)=L(\tau) x_t + \varepsilon f(t, x_t),
$$
where $L: \mathbb{R}_+\rightarrow \mathcal{L}(C; \mathbb{R})$ and $f: \mathbb{R} \times C \rightarrow \mathbb{R}$ are, respectively, linear and nonlinear operators, and $\varepsilon>0$ is a small enough parameter. The theoretical results are applied to a Parkinson's disease model, where the obtained conclusions are illustrated by numerical simulations. 
\end{abstract}

\maketitle
\section{Introduction}\label{Introduccion}

In the last decades the scalar nonlinear differential equations with several delays have become of interest in many fields of science. In particular, in the context of medical treatment, the neuroscience field has received significant attention  (see, e.g., \cite{Book1,Book2}). Based in time series analysis with nonlinear delay differential equations, Lainscsek \emph{et al.}  noticed that several models in medicine (see, e.g., \cite{Lainscsek1,Lainscsek3,Lainscsek2,Lainscsek4,Lainscsek5,Lainscsek6}) can be investigated by the following general class of differential equations with $n$ delays 
\begin{gather}\label{1.1}
\begin{aligned}
\dot{x}(t) &= a_1 x_{\tau_1}+a_2 x_{\tau_2} +a_3 x_{\tau_3}+ \ldots + a_{i-1}x_{\tau_n}+\ldots  \\
        & \ \ +  a_i x^2_{\tau_1}+a_{i+1} x_{\tau_1}x_{\tau_2}+a_{i+2} x_{\tau_1}x_{\tau_3}+ \ldots   \\
        & \ \  +  a_{j-1}x^2_{\tau_n}+a_j x^3_{\tau_1}+a_{j+1}x^2_{\tau_1}x_{\tau_2}+ \ldots \\
        & \ \ \ \vdots \\
        &  \ \ +  \ldots a_l x^m_{\tau_n}. 
\end{aligned}
\end{gather}
Here $a_i \in \mathbb{R}$ and $x_{\tau_i}=x(t-\tau_i)$ where $\tau_i\geq 0 \  (i=1, \ldots , n)$. In particular, Lainscsek \emph{et al.} \cite{Lainscsek2} proposed an algorithm to find a scalar delay differential equation for 
the classification of the finger tapping movement in Parkinson's disease (PD) patients and, as a consequence,  the following equation with 19 terms was considered as a general class of models 
\begin{gather}\label{1.2}
\begin{aligned}
\dot{x}(t) &= a_1 x_{\tau_1}+a_2x_{\tau_2} +a_3 x_{\tau_3}+ a_4 x^2_{\tau_1}+ a_5 x_{\tau_1}x_{\tau_2} \\
        & +a_6x_{\tau_1}x_{\tau_3} + a_7x^2_{\tau_2} +a_8 x_{\tau_2}x_{\tau_3}+a_9x^2_{\tau_3} + a_{10}x^3_{\tau_1} \\
        & +a_{11}x^2_{\tau_1}x_{\tau_2}+a_{12}x^2_{\tau_1}x_{\tau_3} +a_{13}x_{\tau_1}x^2_{\tau_2} \\
        & + a_{14}x_{\tau_1}x_{\tau_2}x_{\tau_3}+a_{15}x_{\tau_1}x^2_{\tau_3}+a_{16}x^3_{\tau_2} \\
        & + a_{17} x^2_{\tau_2}x_{\tau_3}+a_{18}x_{\tau_2}x^2_{\tau_3} +a_{19} x^3_{\tau_3}. 
\end{aligned}
\end{gather}
The above model is useful for assessing the progression of the disease and, in particular, the coefficients in \eqref{1.2} can be used to asses the severity of the disease. Although the above-mentioned models are very important in medical treatment, to the best of our knowledge, delayed models with several delays have not been sufficiently investigated yet. Depending on the data type  \cite{Lainscsek4}, some of the coefficients of the model \eqref{1.2} could be zero or be considered as small parameters, therefore, the results obtained in this paper are applied to find the sufficient conditions for the existence of periodic solutions to a particular case of \eqref{1.2} with a small parameter. 

Averaging is a classical method in the analysis of nonlinear differential equations with small parameters. For ordinary differential equations, many theorems deals with the following initial value problem
\begin{equation}\label{1.3}
\dot{x}(t)=\varepsilon F(t,x) + \varepsilon^2 R(t,x,\varepsilon), \qquad x(0)=x_0,
\end{equation}
where $\varepsilon>0$ is a small parameter and $F, R$ are $T$-periodic in the first argument. According to periodic averaging theorems (see, e.g., \cite{Hale1,Verhulst, Verhulst2}) many dynamical information of the problem \eqref{1.3} can be obtained from the analysis of the following autonomous differential equation
\begin{equation}\label{1.4}
\dot{y}(t)= \frac{1}{T} \int_0^T F(t, y)dt, \qquad y(0)=x_0,
\end{equation}
and, in particular, the existence of $T$-periodic solutions to \eqref{1.3} can be investigated from the analysis of equilibria of the equation \eqref{1.4}. For more general settings, for instance, in \cite{Mesquita2,Lunel,Lakrib2,Lakrib,Mesquita,Halanay,Halesmallparameter} the authors studied the following functional differential equation 
\begin{equation}\label{1.5}
\dot{x}(t)=\varepsilon f(t,x_t),
\end{equation}
where $x_t(\theta)=x(t+\theta)$ for $\theta\in[-r,0]$ and, generally speaking, $ f: \mathbb{R}\times C \rightarrow  \mathbb{R}$ is a nonlinear function. As for ordinary differential equations, good approximations of solutions to \eqref{1.5} can be obtained from solutions of an averaged autonomous functional differential equation. 

In this paper, based on some of the ideas and results in \cite{Faria, Hale2, Hale1,Hsmith, Halesmallparameter,Verhulst}, we study the existence and stability of periodic solutions to the following functional differential equation 
\begin{equation}\label{1.6}
\dot{x}(t)= L(\tau)x_t + \varepsilon f(t, x_t)
\end{equation}
when $\varepsilon$ varies. Here $L: \mathbb{R}_+ \rightarrow \mathcal{L}(C; \mathbb{R})$ is a scalar linear functional, $f: \mathbb{R} \times C \rightarrow \mathbb{R}$ is a nonlinear operator, nonautonomous, $T$-periodic in the first argument and has a continuous second Frech\'et derivative with respect to the second argument and $\varepsilon>0$ is a small enough parameter. First, in the case that the nonperturbed equation
\begin{equation}\label{1.7}
\dot{x}(t)=L(\tau)x_t,
\end{equation}
undergoes a Hopf bifurcation at the critical value $\tau=\tau_0$, we associate to the equation \eqref{1.6} an abstract ordinary differential equation of the form \eqref{1.3}. Secondly, by using classical averaging techniques for ordinary differential equations (see, e.g.,  \cite{Verhulst}), we state sufficient conditions on the parameters to stablish the existence of periodic solutions  for the abstract ordinary differential equation with a small enough parameter $\varepsilon>0$. The stability of such solutions is investigated, as well. The above procedure allow us to describe adequately the dynamic of \eqref{1.6} in a finite dimensional space by restricting the spectrum of the infinitesimal generator of the nonperturbed equation \eqref{1.7}. 

As a consequence of our theoretical results, we investigate a delayed model with a small parameter motivated by the Parkinson's disease model \eqref{1.2} where, in particular, the linear operator $L$ takes the form $L(\tau_2)x_t=-a_1x(t-\tau_1)-a_2x(t-\tau_2)$. 
  
The paper is structured as follows. In Section \ref{s2}, in the case that the equation \eqref{1.7} undergoes a Hopf bifurcation at some critical parameter, we present the construction of an abstract ordinary differential equation that describes the dynamics of \eqref{1.6} and, by applying the classical averaging theory for ordinary differential equations, we prove the existence of periodic solutions for the equation \eqref{1.6} when $\varepsilon$ is small enough. In Section \ref{s3}, by applying the results of Section \ref{s2}, we prove the existence of periodic solutions for a Parkinson's disease model with a small parameter motivated by \eqref{1.2}. To verify the existence of the periodic solution, we present numerical simulations using the Matlab ddesd Package. The last section is devoted to a discussion of the results. 

\section{Abstract ordinary differential equation and averaging theory}\label{s2}

In this section, by applying the general theory in \cite{Faria, Hale2, Hale1, Hsmith, Halesmallparameter}, we introduce a convenient coordinate system to obtain an abstract ordinary differential equation which is equivalent to the general equation \eqref{1.6} in the case that \eqref{1.7} undergoes a Hopf bifurcation. Then, by applying the averaging theory for ordinary differential equations in \cite{Verhulst}, we investigate the existence and stability of periodic solutions for such a ordinary differential equation with a small enough parameter $\varepsilon>0$ and, as a  consequence, we obtain the main result about the existence and stability of periodic solutions for the equation \eqref{1.6}. 

Let $r>0$ and define the phase space $C\stackrel{def}{=} C([-r,0]; \mathbb{R})$ equipped with the sup norm, and consider the following Banach space 
\begin{equation*}
BC\stackrel{def}{=} \lbrace \phi:[-r, 0] \rightarrow \mathbb{R}: \phi \ \textit{is continuous on } [-r, 0[, \exists \lim_{\theta \rightarrow 0^-} \phi (\theta) \in \mathbb{R}\rbrace
\end{equation*}
to define the scalar linear functional $L: \mathbb{R}_+ \rightarrow \mathcal{L}(C; \mathbb{R})$ on $BC$ as
\begin{equation*}
L(\tau)\phi \stackrel{def}{=} \int_{-\tau}^0 \phi (\theta) d\eta(\theta),
\end{equation*}
where $\eta$ is a bounded variation function. Let the linear autonomous equation
\begin{equation}\label{2.1}
\dot{x}(t)=L(\tau)x_t
\end{equation}
and consider the scalar perturbed equation with a small parameter $\varepsilon>0$
\begin{equation}\label{2.2}
\dot{x}(t)=L(\tau)x_t+\varepsilon f(t,x_t).
\end{equation}
We now formulate some of the assumptions on $f$ to state the main results.
\begin{enumerate}
\item[(H.1)] $f$ transforms $ D \subset \mathbb{R}\times C $ into $ \mathbb{R}$ and is a $T$-periodic function in $t$ uniformly with respect to $\phi$ (in some subset of $C$).
\item[(H.2)] $f$ is a continuous and uniformly bounded function in $(t, \phi) \in D$ and has a continuous second Frech\'et derivative with respect to $\phi$. 
\end{enumerate}
Since our purpose is to study the equation \eqref{2.2} in the case that \eqref{2.1} undergoes a Hopf bifurcation at some critical value $\tau=\tau_0$, for the characteristic equation
\begin{equation}\label{2.3}
h(\lambda,\tau) \stackrel{def}{=} \lambda-L(\tau)e^{\lambda}=0
\end{equation}
we shall assume the following hypotheses:
\begin{enumerate}
\item[(H.3)] There exists a pair of simple characteristics roots $\mu(\tau)\pm i\omega(\tau)$ of \eqref{2.3} such that
\begin{equation*}
\mu(\tau_0)=0, \ \ \omega (\tau_0)  \stackrel{def}{=} \omega^*>0, \ \ \mu'(\tau_0)\neq 0.
\end{equation*}
\item[(H.4)] The characteristic equation $h(\lambda, \tau_0)=0$ has no other roots with zero real parts. 
\end{enumerate}

In view of the above-mentioned hypotheses, throughout this section, the following equation will be considered
\begin{equation}\label{2.4}
\dot{x}(t)=L(\tau_0)x_t+\varepsilon f(t,x_t).
\end{equation}
Let $C^*\stackrel{def}{=}C([0, r]; \mathbb{R})$ and, together with the equation \eqref{2.1},  consider the adjoint equation 
\begin{equation*}
\dot{v}(s)=-\int_{-\tau}^0 v(s-\theta)d\eta(\theta)
\end{equation*}
with respect to the bilinear form $(\cdot, \cdot)$ in $C^* \times C$ defined by 
\begin{equation}\label{2.5}
(\psi, \varphi)= \psi(0)\varphi(0)-\int _{-\tau}^0 \int_0^{\theta} \psi(\xi-\theta)d\eta (\theta) \varphi(\xi)d\xi,  \qquad (\psi, \varphi)  \in C^*\times C. 
\end{equation}
Let $L_0=L(\tau_0), ~ \Lambda=\lbrace i\omega^*, -i\omega^* \rbrace$ and $P$ be the center space of the equation  
\begin{equation}\label{2.6}
\dot{y}(t)=L_0y_t.
\end{equation}
Decomposing $C$ by $\Lambda$ as $C=P\oplus Q$, in contrast to the normal forms theory for delay differential equations (see, e.g, \cite[Chapt. 8]{Hale2}) where is useful to consider complex coordinates, since we are interested in real solutions to \eqref{2.2}, along this paper we will choose real bases $\Phi, \Psi$ for $P$ and his dual $P^*$ respectively as

\begin{eqnarray}\label{2.7}
\begin{array}{lcl}
P = \mbox{span} \Phi,  && \Phi(\theta)= (\sin(\omega^*\theta),~ \cos(\omega^*\theta)), \ \ - \tau_0 \leq \theta \leq 0,  \\
 P^*=\mbox{span} \Psi, && \Psi(s) = \textit{col} \ (\alpha_1 \sin(\omega^*s)+\beta_1\cos(\omega^*s),~\alpha_2\sin(\omega^*s)+\beta_2\cos(\omega^*s)), \\ && \ \ \ \ \ \ \ \ \ \ 0\leq s \leq \tau_0
\end{array} 
\end{eqnarray}
where, in order to simplify the transformations, the parameters $(\alpha_1, \beta_1, \alpha_2, \beta_2)\in \mathbb{R}^4$ which could depend on $(\omega^*, \tau_0)$, can be chosen such that $(\Psi, \Phi)=I$ (the identity matrix). Further, the basis $\Phi$ defined as earlier, satisfies
\begin{equation}\label{2.8}
\Phi(\theta)=\Phi(0)e^{B\theta}, \qquad B=  \left( \begin{array}{cc}
0 & -\omega^* \\ 
\omega^* & 0
\end{array} \right).
\end{equation}
We are in position to introduce a coordinate system in $C$ which allow us to obtain the abstract ordinary differential equation associated to \eqref{2.4} (see, e.g.,   \cite{Faria, Hale2,Hale1}). In fact, for the set of eigenvalues $\Lambda  = \lbrace i \omega^*, -i \omega^* \rbrace$, it is possible decompose the elements of $C$ in a unique way as
\begin{equation}\label{2.9}
\phi = \phi^P+\phi^Q, \qquad (\phi^P, \phi^Q)\in P\times Q
\end{equation}
where $Q=\lbrace \phi \in C: (\Psi, \phi)=0  \rbrace$ and $\phi^P=\Phi b, ~b=\textit{col}~(b_1,b_2)=(\Psi, \phi)$ where 
\begin{equation*}
b_i= \beta_i \phi(0)-\int_{-\tau_0}^0\int_{0}^{\theta} \left[ \alpha_i \sin(\omega^*(\xi-\theta))+ \beta_i \cos(\omega^*(\xi-\theta)) \right] d\eta(\theta)\phi(\xi)d\xi  \qquad ( i=1,2).
\end{equation*}
In the sequel, we consider the equation \eqref{2.4} subject to the following non-negative initial condition and positive initial value
\begin{equation}\label{2.10}
x(t)=\varphi(t), \ \ \varphi(t)\geq 0, \ \ -\tau_0 \leq t <0, \ \ x(0)=x_0>0.
\end{equation}
Let $T(t), \ t\geq 0$ be the semigroup on $C$ generated by the solutions of \eqref{2.6},
then any solution of the problem \eqref{2.2}, \eqref{2.10} satisfies the so-called \emph{variation of constant formula} 
\begin{equation}\label{2.11}
x_t=T(t)\varphi + \varepsilon \int_{0}^t T(t-s)X_0f(s,x_s)ds,
\end{equation} 
where 
\begin{equation*}
X_0(\theta)= \left\lbrace \begin{array}{lcc}
             1, &  &  \theta =0 \\
             \\ 0, &  & -\tau_0\leq \theta < 0
             \end{array}
   \right. 
\end{equation*}
Further, in view of \eqref{2.11} and \eqref{2.9}, we have that $x_t=\Phi y(t)+x_t^Q$ satisfies 
\begin{eqnarray}\label{2.12}
\left\lbrace
\begin{array}{lcl}
\dot{y}(t) &=& By(t)+ \varepsilon \Psi (0) f(t, \Phi y(t)+ x_t^Q), \ \  y(0)=(\Psi, \varphi), \\
x_t^Q &=& \displaystyle{T(t)\varphi ^Q+\varepsilon  \int_0^t T(t-s)X_0^Q f(s, \Phi y(s)+x_s^Q) ds}
\end{array}
\right.
\end{eqnarray}
where $B$ is defined in \eqref{2.8}, $~ y(t)=(\Psi, x_t), ~ X_0=X_0^P+X_0^Q$ where $~ X_0^P=\Psi(0)$, and the initial condition satisfies $\varphi=\Phi(\Psi, \varphi)+\varphi^Q$. The integral in \eqref{2.12} is interpreted as a regular integral for $x_t^Q(\theta)$ and each $\theta \in [-\tau_0, 0]$.

\begin{remark}\label{Remark2.1}
It is well-known that for any $\varphi^Q \in Q$, there exists positive constants $K, \alpha$ such that 
\begin{equation*}
|| T(t)\varphi ^Q || \leq Ke^{-\alpha t}|| \varphi ^Q||, \qquad t \geq 0. 
\end{equation*}
Thus, any bounded solution of system \eqref{2.12} satisfies $x_t^Q=O(\varepsilon)$ as $\varepsilon \rightarrow 0$ and hence, since we are interested in $\varepsilon$ sufficiently small, we study only the terms of order $\varepsilon$. 
\end{remark}

According to Remark \ref{Remark2.1}, now the investigation is focused on the initial value problem
\begin{equation}\label{2.13}
\dot{y}(t)=By(t)+\varepsilon\Psi (0)f(t, \Phi y(t)), \qquad y(0)=(\Psi, \varphi).
\end{equation}
In order to study the existence of periodic orbits to \eqref{2.13} using the averaging theory, in the sequel, we proceeds with the introduction of convenient polar coordinates and the application of averaging techniques as Taylor expansions around $\varepsilon=0$. As a first step, we present the \emph{first order theorem} from the averaging theory which we need to prove our principal results (see, e.g., \cite[Thm. 11.5-11.6, pp. 158-159]{Verhulst}).

Consider the ordinary differential system
\begin{equation}\label{2.14}
\dot{x}(t)=\varepsilon F(t,x)+ \varepsilon ^2 R (t,x,\varepsilon),
\end{equation}
where $\displaystyle{F: \mathbb{R}\times D \rightarrow \mathbb{R}^n,
R:  \mathbb{R}\times D \times [0, \varepsilon_{0} ]
\rightarrow \mathbb{R}^{n}}$ with $t \geq 0$ and $D$ is an open
subset of $\mathbb{R}^{n}$. Suppose that
\begin{enumerate}
\item[(a)] The vector functions $\displaystyle{F, R, D_{x}F, D^{2}_{x}F, D_{x}R}$ are continuous and bounded by a constant $M$ (independent of $\varepsilon$) in $[0, \infty [ \times D$, $ \varepsilon \in [0, \varepsilon_{0}]$.
\item[(b)] $F$ and $R$ are $T$-periodic functions in $t$ ( $T$ independent of $\varepsilon$).
\end{enumerate}
The averaged system associated to the system \eqref{2.14} is defined
by
\begin{equation}\label{2.15}
\dot{z}(t)= \varepsilon F^0(z),
\end{equation}
where 
\begin{equation}\label{2.16}
 F^0(z)= \frac{1}{T} \int_{0}^{T} F (s,z)ds.
\end{equation}
\begin{teo}[cf. \cite{Verhulst}]\label{teorema2} Assume that the previous hypotheses $(a)$ and $(b)$ hold and also
that the equation \eqref{2.15} has an equilibrium solution $z^*$
such that $detD_{z}F^0(z^*) \neq 0$, where $D_{z}F^0 (\cdot)$ is a
Jacobian matrix of $F^0$. Then
\begin{enumerate}
\item[(i)] For $\varepsilon$ sufficiently small, there exists a $T$-periodic solution $g(t, \varepsilon) $ of the system \eqref{2.14} such that
$\displaystyle{\lim_{\varepsilon \rightarrow 0} g (t,
\varepsilon)=z^*}$.

\item[(ii)] If the eigenvalues of the critical point $z=z^*$ of the equation \eqref{2.15} all have negative real parts, the corresponding periodic solution $g(t, \varepsilon)$ of \eqref{2.14} is asymptotically stable for $\varepsilon$ sufficiently small. If one of the eigenvalues has positive real part, $g(t, \varepsilon)$ is unstable. 
\end{enumerate}
\end{teo}

From \eqref{2.13}, we will study the equation \eqref{2.4} through the analysis of the  initial value problem
\begin{eqnarray}\label{2.17}
\left\lbrace
\begin{array}{lcl}
\dot{y}_1(t)&=& - \omega^* y_2(t)+ \varepsilon \beta_1 f(t, \Phi y(t)), \ \ \ y_1(0)=y_{10}, \\
\dot{y}_2(t)&=& \ \   \omega^* y_1(t)+ \varepsilon \beta_2 f(t, \Phi y(t)), \ \ \ y_2(0)=y_{20}.
\end{array}
\right.
\end{eqnarray}
Here $\Phi y(t)=y_1(t) \sin(\omega^* \theta)+y_2(t)\cos(\omega^* \theta), ~ -\tau_0\leq \theta \leq 0$ and the initial conditions satisfy 
\begin{equation}\label{2.18}
y_{i0}= \beta_i \varphi(0)-\int_{-\tau_0}^0\int_{0}^{\theta} \left[ \alpha_i \sin(\omega^*(\xi-\theta))+ \beta_i \cos(\omega^*(\xi-\theta)) \right] d\eta(\theta)\varphi(\xi)d\xi  \qquad (i=1,2).
\end{equation}
We introduce the polar coordinates $(\rho, \xi) \in \mathbb{R}^+\times \mathbb{S}^1$  by the change of variables
\begin{equation}\label{2.19}
y_1(t)= \rho(t) \sin(\omega^* \xi(t)), \ \ \ \ y_2(t)=-\rho(t) \cos(\omega^* \xi(t)),
\end{equation}
whence \eqref{2.17} becomes
\begin{eqnarray}\label{2.20}
\left\lbrace
\begin{array}{lcl}
\dot{\rho}(t)&=& \varepsilon f(t,\Phi y(t)) (\beta_1 \sin(\omega^* \xi(t))-\beta_2 \cos(\omega^* \xi(t))), \qquad \qquad \rho(0)=\rho_0, \\
\dot{\xi}(t) &=& \displaystyle{1+ \varepsilon f(t, \Phi y(t)) \left(\frac{\beta_1 \cos(\omega^* \xi(t))+\beta_2 \sin(\omega^* \xi(t))}{ \omega^* \rho(t)}\right)}, \ \ \ \xi(0)=\xi_0,
\end{array}
\right. 
\end{eqnarray}
where $\Phi y(t)= \rho(t) (\sin(\omega^* \xi(t)) \sin(\omega^* \theta)-\cos(\omega^* \xi(t))\cos(\omega^* \theta)), ~ -\tau_0\leq \theta \leq 0$, and the above initial conditions satisfy
\begin{equation*}
y_{10}=\rho_0\sin(\omega^*\xi_0), \ \ y_{20}=-\rho_0 \cos(\omega^*\xi_0),
\end{equation*}
where $y_{10}$ and $y_{20}$ are defined in \eqref{2.18}.  
\begin{remark}
Note that the derivatives of \eqref{2.20} are with respect to time $t$ but the system is not periodic in $t$. If instead of $t$, we take the variable $\xi$ as the new independent variable of the system, we obtain the periodicity necessary to apply Theorem \ref{teorema2}. 
\end{remark}
For $\varepsilon >0$ sufficiently small and in a neighborhood of $\rho=0$ we have that $\dot{\xi}\neq 0$, in such a neighborhood we take $\xi$ as the new independent variable and we denote by a prime the derivative with respect to $\xi$. At the same time, expanding in Taylor series around $\varepsilon=0$ we obtain the differential equation 
\begin{equation}\label{2.21}
\rho '(\xi) = \varepsilon F(\xi, \rho(\xi)) + O(\varepsilon^2), \qquad \rho(0)=\rho_0, 
\end{equation}
where 
\begin{eqnarray*}
\begin{array}{lcl}
&&F(\xi, \rho(\xi)) \stackrel{def}{=}  (\beta_1\sin(\omega^* \xi)-\beta_2 \cos(\omega^* \xi))~ f(\xi, \rho (\xi) (\sin(\omega^* \xi)\sin(\omega^* \theta)-\cos(\omega^* \xi)\cos(\omega^* \theta))),  \\
 &&  \qquad \qquad \qquad -\tau_0 \leq \theta \leq 0.
\end{array} 
\end{eqnarray*}
In what follows, we assume that $f$ is a $T=2\pi/\omega^*$ - periodic function in $\xi$. Then, according to \eqref{2.15}-\eqref{2.16}, the averaged equation associated to \eqref{2.21} is 
\begin{equation}\label{2.22}
\rho '(\xi) = \varepsilon F^0(\rho),
\end{equation}
where
\begin{equation*}
F^0(\rho)  \stackrel{def}{=} \frac{~ \omega^*}{2\pi} \int_{0}^{2\pi / \omega^*} F(\xi, \rho)d\xi .
\end{equation*}
\begin{remark}
If $f$ satisfies the above-mentioned conditions (H.1)-(H.2) for $T=2\pi/\omega^*$, then the assumptions of Theorem \ref{teorema2} are fulfilled for the equation \eqref{2.21}. On the other hand, it is easy to check that $G: [0, +\infty[ \times [0, 2\pi[ \rightarrow \mathbb{R}^2$ defined by $G(\rho, \xi)= (\rho \sin(\omega^*\xi), -\rho \cos(\omega^*\xi))$  is a $C^1$-diffeomorphism and, since the system \eqref{2.17} is a special case of \eqref{2.20}, for a small parameter $\varepsilon>0$ the study about the existence and stability of $T$-periodic solutions for the ordinary differential equation \eqref{2.21} will give us equivalent information about $T$-periodic solutions to \eqref{2.4}, \eqref{2.10}.
\end{remark}
According to this last remark and, as a consequence of Theorem \ref{teorema2}, we now state the principal result of this section.
\begin{teo}\label{promediorho}
If there exists an equilibrium $\rho^*>0$ of the equation \eqref{2.22} such that $F^0(\rho^*)=0$ and $D_{\rho}F^0(\rho^*) \neq 0$ then, for $\varepsilon$ sufficiently small, there exists a $~ 2\pi/\omega^*$- periodic solution $x(t)=g(t,\varepsilon)$ of the problem \eqref{2.4},  \eqref{2.10} such that $g(t,0)=\rho^*$. Moreover, such a periodic solution is asymptotically stable (resp. unstable) on the center manifold if the positive equilibrium $\rho^*$ is stable (resp. unstable).
\end{teo}

\section{A two delayed equation and applications}\label{s3}
In this section, first we show two of the main results presented by Piotrowska \cite{Piotrowska} about the existence of Hopf bifurcations in the following particular case of the equation \eqref{2.1} 
\begin{equation}\label{3.1}
\dot{x}(t)=-a_1x(t-\tau_1)-a_2x(t-\tau_2).
\end{equation}
Here $a_1,a_2\in \mathbb{R}$, and the nonnegative delays are independent of each other. Secondly, we apply the above-mentioned results and the results of Section \ref{s2} to prove the existence of periodic solutions to a particular case of the Parkinson's disease model \eqref{1.2}. The results are illustrated by numerical simulations.

Li  \emph{et al.} \cite{LiRuan} investigated the stability of the zero equilibrium of \eqref{3.1} and the existence of Hopf bifurcations considering $\tau_2$ as a bifurcation parameter, where suitable sets of $a_1, a_2$ and $\tau_1$ were stablished to prove their main results. Piotrowska \cite{Piotrowska} presented some remarks and improvement of the results in \cite{LiRuan}. Furthermore, some aditional cases than the cases presented in \cite{LiRuan} were treated in \cite{Piotrowska}. In the sequel, we keep in mind that $\tau_2$ will be the bifurcation parameter. If $a_1a_2=0, \tau_1\tau_2=0$ or $\tau_1=\tau_2$, then the equation \eqref{3.1} is equal to an equation with a single delay, therefore, no interest is posed in this latter cases. 

The characteristic equation associated with \eqref{3.1} is 
\begin{equation}\label{3.2}
h(z)=z+a_1 e^{-z\tau_1}+a_2e^{-z\tau_2}=0.
\end{equation}
By assuming $a_1\neq 0, a_2>0, \lambda=\dfrac{z}{a_1}, \tau_1=\dfrac{r_1}{|a_1|}, \tau_2= \dfrac{r_2}{|a_1|}$ and $a=\dfrac{a_2}{|a_1|},$ then for $a_1>0$ the characteristic equation becomes 
\begin{equation}\label{3.3}
\lambda=-e^{-\lambda r_1}-ae^{-\lambda r_2},
\end{equation}
and for $a_1<0$ becomes
\begin{equation}\label{3.4}
\lambda =e ^{-\lambda r_1} -a e^{-\lambda r_2}.
\end{equation}
If $a_1+a_2=0,$ then $z=0$ is a solution of \eqref{3.2} for all nonnegative delays $\tau_1, \tau_2 $ and, therefore, we do not consider that case. Also, from \cite[Prop. 1]{HaleHuang} if $a_1+a_2<0,$ then the zero solution of \eqref{3.2} is unstable for all delays $\tau_1, \tau_2$. Thus, throughout this section we assume $a_1+a_2>0$. It is easy to check that if $\tau_1=\tau_2=0$ and $a_1+a_2>0$ are fulfilled, then the zero solution of \eqref{3.2} is asymptotically stable and, by continuity, it will be asymptotically stable for sufficiently small delays $\tau_1, \tau_2>0$. 

We are interested in the existence of purely imaginary roots for \eqref{3.3} ($a_1>0$) and \eqref{3.4} ($a_1<0$), hence first we do $\lambda=i\omega (\omega>0)$ in both equations. Secondly, by separating real and imaginary parts, we have that the equations \eqref{3.3} and \eqref{3.4} have purely imaginary roots if the following systems
\begin{eqnarray}\label{3.5}
\begin{array}{lcl}
&&  \cos(\omega r_1)=-a\cos(\omega r_2), \\
 & & \omega - \sin (\omega r_1)= a\sin (\omega r_2),  
\end{array}
\end{eqnarray}
and 
\begin{eqnarray}\label{3.6}
\begin{array}{lcl}
&& \cos (\omega r_1)= a\cos(\omega r_2),  \\
&& \omega + \sin (\omega r_2)= a\sin (\omega r_2),
\end{array}
\end{eqnarray}
are satisfied, respectively. Adding up to square both sides of \eqref{3.5} and \eqref{3.6} we obtain
\begin{equation}\label{3.7}
\sin (\omega r_1) = \dfrac{\omega^2+1-a^2}{2\omega},
\end{equation}
and 
\begin{equation}\label{3.8}
\sin (\omega r_1) = \dfrac{a^2-\omega^2-1}{2\omega},
\end{equation}
respectively. For $\omega \in  ]0, + \infty[$ the geometrical properties of the functions on the right side of \eqref{3.7} and \eqref{3.8} were considered in \cite{Piotrowska} to investigate the existence of simple purely imaginary roots to \eqref{3.2}. In what follows, we present the critical values of the delay $r_2$ to undergoes a Hopf bifurcation around the zero solution of \eqref{3.1}. We refer to \cite{Piotrowska} to the details and notation involved. 

For any $r_1>0$ the equations \eqref{3.7} and \eqref{3.8} have a finite number $m\in \mathbb{N}$ of solutions $\omega_k$ ($k=1, \ldots, m$), then for all $r_1>0$ fixed and for each $\omega_k$ there is a infinite number of delays $r_2>0$ such that
 \begin{equation}\label{3.9}
 \cos(\omega_k r_1) =- a\cos(\omega_kr_2)
 \end{equation}
 and 
\begin{equation}\label{3.10}
\cos(\omega_k r_1)= a \cos (\omega_k r_2),
\end{equation}
respectively. Let us consider
\begin{equation*}
r_2^{k,+}=\min \lbrace r_2\in \mathbb{R}_+: \eqref{3.9} \  \mbox{is fulfilled} \rbrace
\end{equation*}
and 
\begin{equation*}
r_2^{k,-}= \min \lbrace r_2\in \mathbb{R}_+ : \eqref{3.10} \  \mbox{is fullfiled} \rbrace ,
\end{equation*}
to define 
\begin{equation}\label{3.11}
r_2^{0,+} = \min \lbrace r_2^{k,+}: k=1, \ldots , m \rbrace
\end{equation}
and 
\begin{equation}\label{3.12}
r_2^{0,-}=\min \lbrace r_2^{k,-} : k=1, \ldots, m \rbrace .
\end{equation}
In both of the above-mentioned cases ($a_1>0$ and $a_1<0$), we have that for some $k=1, \ldots ,m$ there exist $r_2^{k, \pm }=r_2^{0, \pm}$, then for that $k$ we put $\omega ^*=\omega_k$. On the other hand, define
\begin{equation*}
\Omega  \stackrel{def}{=} \bigcup_{l \in \mathbb{N}\cup \lbrace 0 \rbrace} [2l\pi , \frac{\pi}{2}+2l\pi].
\end{equation*}
Now we present two of the most important results in \cite{Piotrowska} about the existence of Hopf bifurcations to \eqref{3.1} for $a_1>0$ and $a_1<0$ respectively. 

\begin{teo}[cf. \cite{Piotrowska}]\label{Theorem3.1}
Let $0<a_1<a_2$ and $\tau_2^{0,+}= \dfrac{r_2^{0,+}}{a_1}$, where $r_2^{0,+}$ is defined by \eqref{3.11}, then
\begin{enumerate}
\item[i)] For $\tau_2 \in [0, \tau_2^{0,+}[$ the trivial solution to \eqref{3.1} is asymptotically stable.

\item[ii)] If $\omega^*\tau_2^{0,+} a_1 \in \Omega$, then the equation \eqref{3.1} undergoes a Hopf bifurcation when $\tau_2=\tau_2^{0,+}$. 
\end{enumerate}
\end{teo}

\begin{teo}[cf. \cite{Piotrowska}]\label{Theorem3.2}
 Let $a_1<0, a_2>|a_1|$ and $\tau_2^{0,-}=\dfrac{r_2^{0,-}}{|a_1|}$, where $r_2^{0,-}$ is defined by \eqref{3.12}, then
\begin{enumerate}

\item[i)] For $\tau_2 \in [0, \tau_2^{0,-}[$ the trivial solution to \eqref{3.1} is asymptotically stable.

\item[ii)] If $\omega^*\tau_2^{0,-}|a_1|\in \Omega$, then the equation \eqref{3.1} undergoes a Hopf bifurcation when $\tau_2=\tau_2^{0,-}$. 
\end{enumerate}
\end{teo}

In the sequel, in view of the notation and conclusions of Section \ref{s2}, we consider that the equation \eqref{3.1} satisfies the hypotheses of Theorems \ref{Theorem3.1} or \ref{Theorem3.2}, i.e., for $\tau_2=\tau_2^0 \in \lbrace \tau_2^{0,+}, \tau_2^{0,-} \rbrace $ the equation \eqref{3.1} satisfies the conditions (H.3) and (H.4) mentioned in Section \ref{s2}. Whitout loss of generality we assume that $\tau_1<\tau_2^0$ to define the linear functional
\begin{equation*}
L(\tau_2)\phi = \int_{-\tau_2^0}^0 \phi(\theta)d\eta (\theta),
\end{equation*}
where 
\begin{eqnarray*}
\eta (\theta) = \left\lbrace
\begin{array}{lcl}
0 &\mbox{if}& \qquad \ \ \ \ \theta=-\tau_2^0, \\
-a_2 &\mbox{if}& \ -\tau_2^0<\theta \leq \tau_1 , \\
 -(a_1+a_2) &\mbox{if}&  \ -\tau_1< \theta \leq 0 .
\end{array}
\right. 
\end{eqnarray*}
Furthermore, the bilinear form \eqref{2.5} is
\begin{equation*}
(\psi, \varphi)=\psi(0)\varphi(0)-\int_{-\tau_2^0}^0 [a_1\psi(\xi+\tau_1)+a_2\psi(\xi+\tau_2^0)]\varphi(\xi)d\xi
\end{equation*}
and the coefficients $\alpha_1, \beta_1, \alpha_2, \beta_2$ presented in \eqref{2.7} are
\begin{gather*}
\begin{aligned}
  \alpha_1 & \stackrel{def}{=} \dfrac{4-2\sin^2(\omega^*\tau_2^0)}{\omega^{*^2}\tau_2^{0^2}+\sin^2(\omega^*\tau_2^0)} , \\
  \beta_1  & \stackrel{def}{=} \dfrac{2\omega^*\tau_2^0+\sin(2\omega^*\tau_2^0)}{\omega^{*^2}\tau_2^{0^2}+\sin^2(\omega^*\tau_2^0)} ,\\
 \alpha_2 &  \stackrel{def}{=} \dfrac{\sin(2\omega^*\tau_2^0)-2\omega^*\tau_2^0}{\omega^{*^2}\tau_2^{0^2}+\sin^2(\omega^*\tau_2^0)} , \\
 \beta_2&\stackrel{def}{=}\dfrac{2\sin^2(\omega^*\tau_2^0)}{\omega^{*^2}\tau_2^{0^2}+\sin^2(\omega^*\tau_2^0)}.
\end{aligned}
\end{gather*}
Such a coefficients are such that $(\Psi, \Phi)=I$, where $\Phi$ and $\Psi$ are the bases for the center space $P$ and his dual $P^*$ respectively. In what follows we apply the results of Section \ref{s2} and the results of \cite{Piotrowska} presented above to investigate a particular case of \eqref{1.1}. 

Most of the terms in \eqref{1.1} could be set to zero depending on the data type \cite{Lainscsek4}, therefore, now we investigate the following slightly modified Parkinson's disease model with negative feedback, three delays and a small parameter
\begin{equation}\label{3.13}
\dot{x}(t)=-a_1x(t-\tau_1)-a_2x(t-\tau_2)+ \varepsilon (a_4x^3(t-\tau_3)-a_3 x (t-\tau_3)).
\end{equation}
Here the real parameters $a_1, a_2, \tau_1$ and  $\tau_2$ satisfy the conditions and conclusions of Theorems \ref{Theorem3.1} or \ref{Theorem3.2}. It is worth mentioning here that, in view of \cite{Lainscsek1}, the coefficients $a_1, a_2, a_3, a_4$ could be positive or negative depending on the data type, and the delays could be choosed to separate PD on or off medication from controls (individuals without the disease), and to separate PD on from PD off medication. 

Let $\tau_2=\tau_2^0 \in \lbrace \tau_2^{0,+}, \tau_2^{0,-} \rbrace$, then in view of \eqref{2.19} we have
\begin{eqnarray*}
x(t-\tau_1) &=& \Phi(-\tau_1)y(t)= -\rho (\sin(\omega^*\tau_1)\sin(\omega^*\xi)+\cos(\omega^*\tau_1)\cos(\omega^*\xi)), \\
x(t-\tau_2^0) &=&  \Phi(-\tau_2^0)y(t)= -\rho (\sin(\omega^*\tau_2^0))\sin(\omega^*\xi)+\cos(\omega^*\tau_2^0)\cos(\omega^*\xi)), \\ 
x(t-\tau_3) &=&  \Phi(-\tau_3)y(t)= -\rho (\sin(\omega^*\tau_3))\sin(\omega^*\xi)+\cos(\omega^*\tau_3)\cos(\omega^*\xi)).
\end{eqnarray*}
Thus, the averaged equation \eqref{2.22} is
\begin{equation*}
\rho'(\xi)=  \frac{1}{8} \rho (3a_4 \rho^2 -4a_3)(\beta_2\cos(\omega^*\tau_3)-\beta_1\sin(\omega^*\tau_3)). 
\end{equation*}
Assume that $a_3 a_4>0$ and the delay $\tau_3>0$ is such that $\beta_2\cos(\omega^*\tau_3)-\beta_1\sin(\omega^*\tau_3)\neq 0$, then the averaged equation has a positive equilibrium $\rho^*=\sqrt{(4a_3)/(3a_4)}$ such that $D_{\rho}F^0(\rho^*)=a_3(\beta_2\cos(\omega^*\tau_3)-\beta_1\sin(\omega^*\tau_3))$. Then if  $a_3(\beta_2\cos(\omega^*\tau_3)-\beta_1\sin(\omega^*\tau_3))<0$ the equilibrium point is asymptotically stable and unstable if $a_3(\beta_2\cos(\omega^*\tau_3)-\beta_1\sin(\omega^*\tau_3))>0$. Consequently, according to Theorem \ref{promediorho}, for $\varepsilon$ sufficiently small if $a_3a_4>0$ we have the existence of a $2\pi/\omega^*$-periodic solution $g(t, \varepsilon)$ to \eqref{3.13} such that $g(t,0)=\rho^*$. Further, such a periodic solution is asymptotically stable (resp. unstable) on the center manifold  if $a_3(\beta_2\cos(\omega^*\tau_3)-\beta_1\sin(\omega^*\tau_3))<0$ (resp. $a_3(\beta_2\cos(\omega^*\tau_3)-\beta_1\sin(\omega^*\tau_3))>0$). 

\begin{figure}[H]
\centering
\subfigure[$\varepsilon=0.1$]{\scalebox{0.24}{\includegraphics{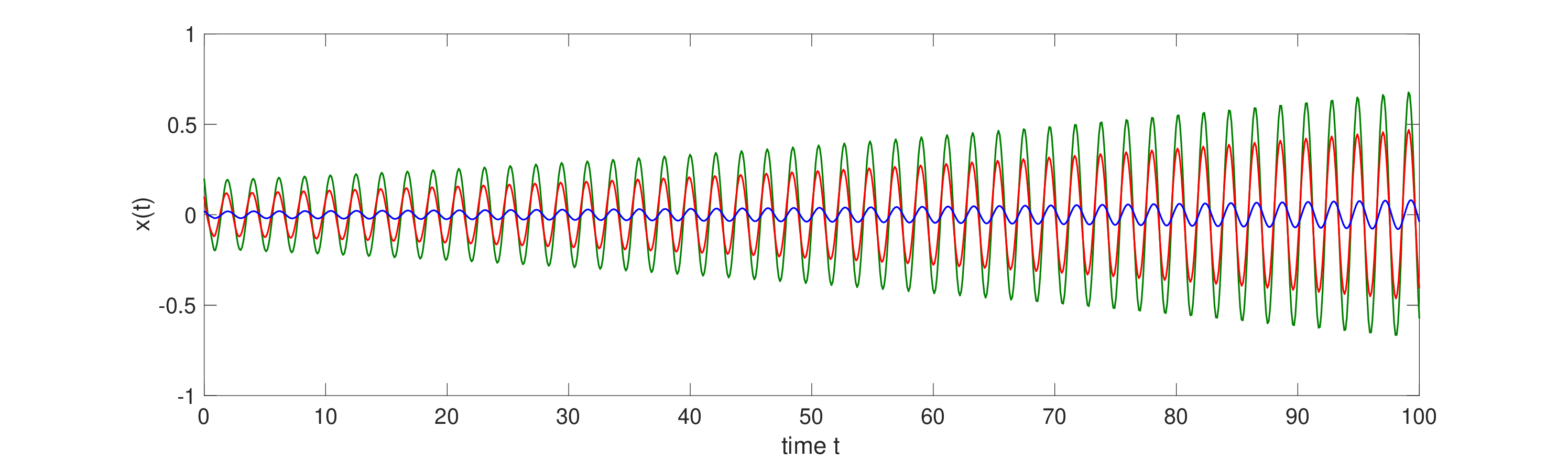}}}
\hspace{0mm}
\subfigure[ $\varepsilon=0.001$]{\scalebox{0.24}{\includegraphics{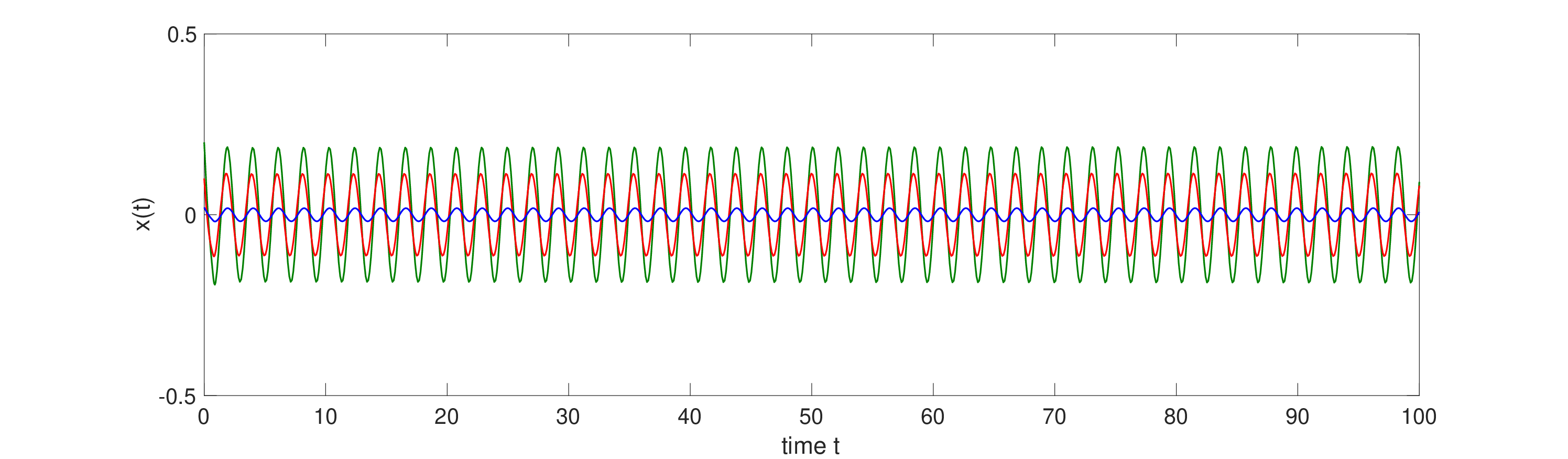}}}
\caption{~ Dynamic behaviour of \eqref{3.14} for $a_3= a_4=1$, where the initial functions are: $\varphi_1(\theta)=0.2 e^{\theta} $ (green, outside), \ $\varphi_2(\theta)=0.05(\cos(\theta)+1)$ (red, middle), $\varphi_3(\theta)=0.02(\sin(\theta)+1)$ (blue, inside) for $\theta\in [-\tau_3,0]$}
\label{F2}
\end{figure}

\begin{figure}[H]
\centering
\subfigure[$\varepsilon=0.1$]{\scalebox{0.24}{\includegraphics{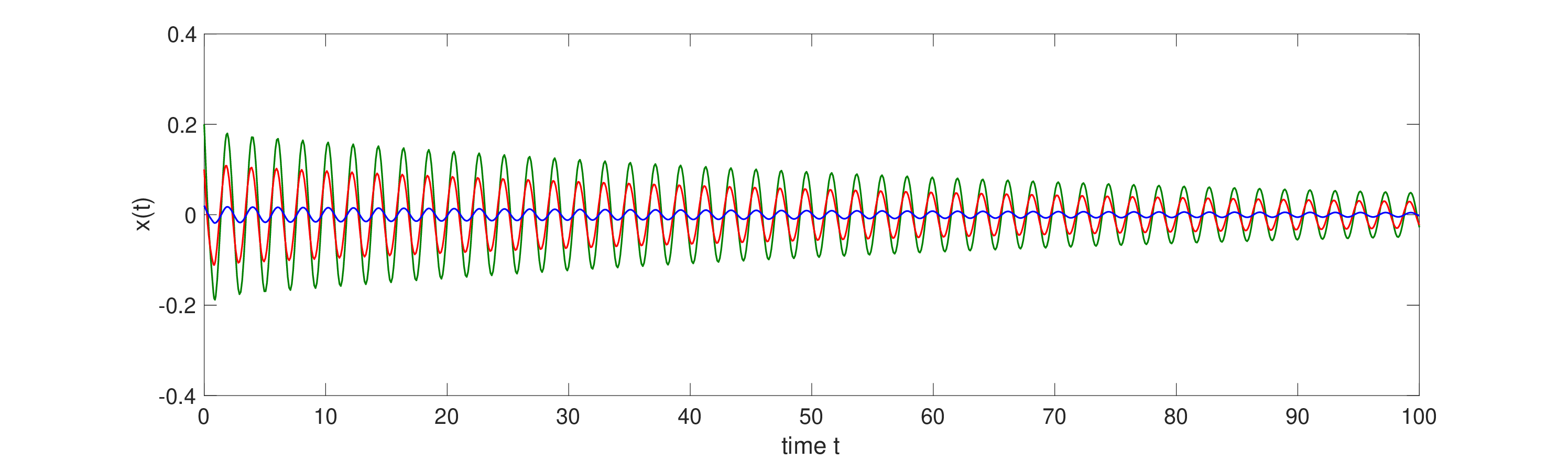}}}
\hspace{0mm}
\subfigure[ $\varepsilon=0.001$]{\scalebox{0.24}{\includegraphics{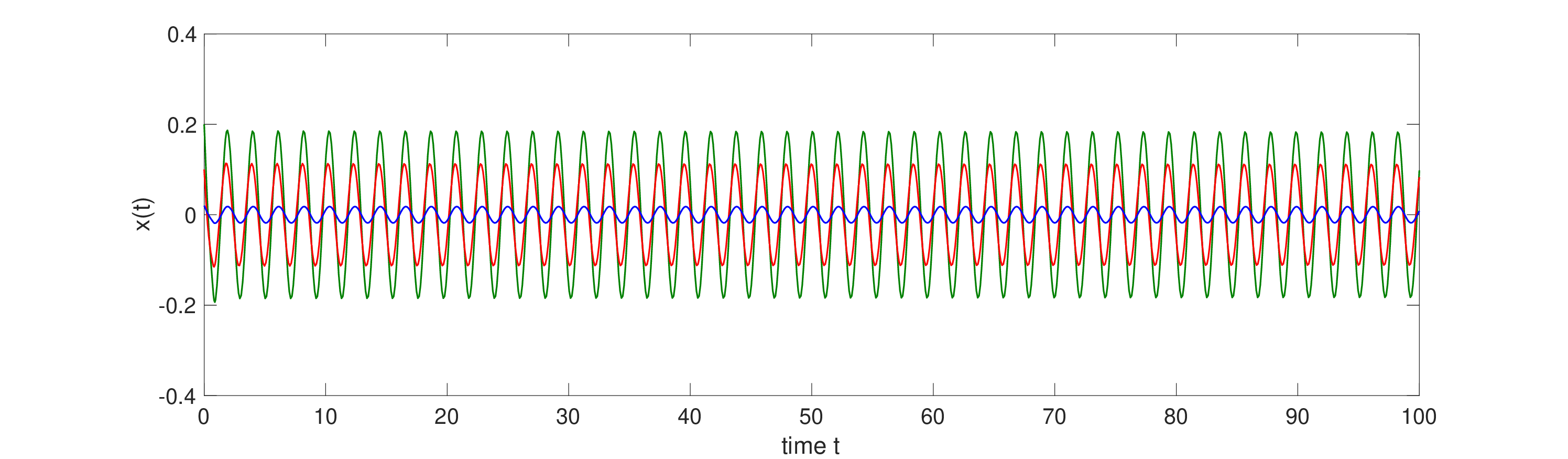}}}
\caption{~ Dynamic behaviour of \eqref{3.14} for $a_3=a_4=-1$, where the initial functions are: $\varphi_1(\theta)=0.2 e^{\theta} $ (green, outside), \ $\varphi_2(\theta)=0.05(\cos(\theta)+1)$ (red, middle), $\varphi_3(\theta)=0.02(\sin(\theta)+1)$ (blue, inside) for $\theta\in [-\tau_3,0]$}
\label{F3} 
\end{figure}

Let us consider a particular case of \eqref{3.13}, namely the equation
\begin{equation}\label{3.14}
\dot{x}(t)=-2x(t-\tau_1)-3x(t-\tau_2)+\varepsilon(x^3(t-\tau_3)-x(t-\tau_3)). 
\end{equation}
As usual, by a solution of the equation \eqref{3.14} we understand an continuously differentiable function $x$ which satisfies the problem \eqref{3.14}, \eqref{2.10} for $t\geq 0$. According to Theorem \ref{promediorho}, for $\tau_2=\tau_2^0$ and $\varepsilon$ sufficiently small, a periodic solution arises around the zero equilibrium. In particular, for $\tau_1=0.113279$ we obtain $\omega^*=3$ and $\tau_2^0=0.750157$. For $\tau_3=1.2$, in Figure \ref{F2} (resp. Figure \ref{F3}) we present a numerical simulation for $a_3=a_4=1$ (resp. $a_3=a_4=-1$), where if $a_3=a_4=\pm 1$ then $D_{\rho}F^0(\rho^*)=\pm 0.08362$, i.e., if $a_3=a_4=1$ (resp. $a_3=-1$) the periodic solution is unstable (resp. asymptotically stable) on the center manifold.

\section{Discussion}
In this paper, we study a general class of scalar delay differential equations with a small parameter. By applying the Hopf bifurcation theorem and the so-called averaging theory for ordinary differential equations, we prove the existence of periodic solutions around the zero solution of perturbed models when the perturbation is sufficiently small. An equation arising from a Parkinson's disease model is investigated to apply our theoretical results, where the numerical simulations illustrate our investigation. Our results can be applied to a wide class of neuroscience models with several delays as long as the conditions of our main results are satisfied. To the best of the author's knowledge, neuroscience models with several delays have not been thoroughly investigated so far, therefore this paper contribute in the investigation of dynamical behaviour of such a models.  

\section*{Acknowledgements}

J. Oyarce acknowledges support from Chilean National Agency for Research an Development (PhD. 2018-21180824).

\end{document}